\magnification=\magstephalf
\def\to{\ \longrightarrow\ }

\def\nl{\hfill\break}

\def\hexnumber#1{\ifcase#1 0\or 1\or 2\or 3\or 4\or 5\or 6\or 7\or 8\or
 9\or A\or B\or C\or D\or E\or F\fi}
%
%
\font\twelvemsa=msam10 scaled 1200   
\font\tenmsa=msam10                  
\font\ninemsa=msam9            \font\sevenmsa=msam7
\font\sixmsa=msam6             \font\fivemsa=msam5
%
%
\newfam\msafam                 \textfont\msafam=\tenmsa
\scriptfont\msafam=\sevenmsa   \scriptscriptfont\msafam=\fivemsa
\edef\hexa{\hexnumber\msafam}        
\def\msa{\fam\msafam\tenmsa}         
%
%
\font\twelvemsb=msbm10 scaled 1200   
\font\tenmsb=msbm10                  
\font\ninemsb=msbm9            \font\sevenmsb=msbm7
\font\sixmsb=msbm6             \font\fivemsb=msbm5
%
\newfam\msbfam                 \textfont\msbfam=\tenmsb       
\scriptfont\msbfam=\sevenmsb   \scriptscriptfont\msbfam=\fivemsb
\edef\hexb{\hexnumber\msbfam}        
\def\msb{\fam\msbfam\tenmsb}         
%
%
\font\twelveeufm=eufm10 scaled 1200  
\font\teneufm=eufm10                 
\font\nineeufm=eufm9           \font\seveneufm=eufm7
\font\sixeufm=eufm6            \font\fiveeufm=eufm5
%
\newfam\eufmfam                \textfont\eufmfam=\teneufm
\scriptfont\eufmfam=\seveneufm \scriptscriptfont\eufmfam=\fiveeufm
\edef\hexf{\hexnumber\eufmfam}      
\def\frak{\fam\eufmfam\teneufm}     
%
%
%
\font\twelverm=cmr10 scaled 1200    
\font\ninerm=cmr9                   
\font\sixrm=cmr6   
%
\font\twelvei=cmmi10 scaled 1200    
\font\ninei=cmmi9                   
\font\sixi=cmmi6  
%
\font\twelvesy=cmsy10 scaled 1200   
\font\ninesy=cmsy9                  
\font\sixsy=cmsy6  
%
\font\twelvebf=cmbx10 scaled 1200   
\font\ninebf=cmbx9                  
\font\sixbf=cmbx6  
%
%
\font\twelveit=cmti10 scaled 1200   
\font\nineit=cmti9                  
%
\font\twelvesl=cmsl10 scaled 1200   
\font\ninesl=cmsl9                  
%
\font\twelvett=cmtt10 scaled 1200   
\font\ninett=cmtt9                  
%
%
%
%
\def\small{%
%
%
\textfont0=\ninerm \scriptfont0=\sixrm \scriptscriptfont0=\fiverm
\def\rm{\fam0\ninerm}        
%
%
\textfont1=\ninei \scriptfont1=\sixi \scriptscriptfont1=\fivei
%
%
\textfont2=\ninesy \scriptfont2=\sixsy \scriptscriptfont2=\fivesy
%
%
\textfont3=\tenex \scriptfont3=\tenex \scriptscriptfont3=\tenex
%
%
\textfont\bffam=\ninebf \scriptfont\bffam=\sixbf
\scriptscriptfont\bffam=\fivebf \def\bf{\fam\bffam\ninebf}%
%
%
\textfont\itfam=\nineit \def\it{\fam\itfam\nineit}%
\textfont\slfam=\ninesl \def\sl{\fam\slfam\ninesl}%
\textfont\ttfam=\ninett \def\tt{\fam\ttfam\ninett}%
%
%
%
\textfont\msafam=\ninemsa \scriptfont\msafam=\sixmsa
\scriptscriptfont\msafam=\fivemsa \def\msa{\fam\msafam\ninemsa}%
%
%
\textfont\msbfam=\ninemsb \scriptfont\msbfam=\sixmsb
\scriptscriptfont\msbfam=\fivemsb \def\msb{\fam\msbfam\ninemsb}%
%
%
\textfont\eufmfam=\nineeufm  \scriptfont\eufmfam=\sixeufm
\scriptscriptfont\eufmfam=\fiveeufm \def\frak{\fam\eufmfam\nineeufm}%
%
%
%
\normalbaselineskip=11pt
\setbox\strutbox=\hbox{\vrule height8pt depth3pt width0pt}%
%
%
\normalbaselines\rm}    
%
%
%
%
\def\large{%
\textfont0=\twelverm \scriptfont0=\ninerm \scriptscriptfont0=\sevenrm
\def\rm{\fam0\twelverm}%
\textfont1=\twelvei \scriptfont1=\ninei \scriptscriptfont1=\seveni
\textfont2=\twelvesy \scriptfont2=\ninesy \scriptscriptfont2=\sevensy
\textfont3=\tenex \scriptfont3=\tenex \scriptscriptfont3=\tenex
\textfont\bffam=\twelvebf \scriptfont\bffam=\ninebf
\scriptscriptfont\bffam=\sevenbf \def\bf{\fam\bffam\twelvebf}%
\textfont\itfam=\twelveit \def\it{\fam\itfam\twelveit}%
\textfont\slfam=\twelvesl \def\sl{\fam\slfam\twelvesl}%
\textfont\ttfam=\twelvett \def\tt{\fam\ttfam\twelvett}%
\textfont\msafam=\twelvemsa \scriptfont\msafam=\ninemsa
\scriptscriptfont\msafam=\sevenmsa \def\msa{\fam\msafam\twelvemsa}         
\textfont\msbfam=\twelvemsb \scriptfont\msbfam=\ninemsb
\scriptscriptfont\msbfam=\sevenmsb \def\msb{\fam\msbfam\twelvemsb}         
\textfont\eufmfam=\twelveeufm  \scriptfont\eufmfam=\nineeufm
\scriptscriptfont\eufmfam=\seveneufm \def\frak{\fam\eufmfam\teneufm}
\normalbaselineskip=15pt
\setbox\strutbox=\hbox{\vrule height11pt depth4pt width0pt}%
\normalbaselines\rm}%
%
\def\Bbb{\msb}

%

%
\mathchardef\plussquare="0\hexa01
\mathchardef\nge="3\hexb0B
\mathchardef\maltesecross="0\hexa7A
\mathchardef\del="0\hexf01
%
\def\ninepoint{\small}

\input epsf
\overfullrule=0pt

\font\Bbb=msbm10

\font\secfont=cmbx10

\font\nam=cmr8
\font\aff=cmti8
\font\refe=cmr9

\mathchardef\square="0\hexa03
\def\qed{\hfill$\square$\par\rm}
\def\np{\vfill\eject}
\def\boxing#1{\ \lower 3.5pt\vbox{\vskip 3.5pt\hrule \hbox{\strut\vrule \ #1 \vrule} \hrule} }

\def\down#1{\ \lower 3.5pt\vbox{\vskip 3.5pt \hbox{\strut \ #1 \vrule} \hrule} }
\def\negdown#1{\ \lower 3.5pt\vbox{\vskip 3.5pt \hbox{\strut  \vrule \ #1 }\hrule} }
\def\adj{\hbox{\rm adj}}

\hsize=6.3 truein
\vsize=9 truein

\baselineskip=13 pt
\parskip=\baselineskip
 1

\parindent=0pt

\def\Z{\hbox{\Bbb Z}}
\def\R{\hbox{\Bbb R}}
\def\C{\hbox{\Bbb C}}
\def\H{\hbox{\Bbb H}}
\def\d{\hbox{\rm d}}
\def\diag{\mathop{\rm diag}\nolimits}







\newif \iftitlepage \titlepagetrue

\def\diagram{\global\advance\diagramnumber by 1
$$\epsfbox{equivalentfig.\number\diagramnumber.eps}$$}
\def\ddiagram{\global\advance\diagramnumber by 1
\epsfbox{equivalentfig.\number\diagramnumber.eps}}

\newcount\diagramnumber
\diagramnumber=0

\newcount\secnum \secnum=0
\newcount\subsecnum \subsecnum=0
\newcount\defnum \defnum=0
\newcount\lemmanum \lemmanum=0
\def\section#1{
                \vskip 10 pt
                \advance\secnum by 1 \subsecnum=0\lemmanum=0
                \leftline{\secfont \the\secnum \quad#1}
                }

\def\subsection#1{
                \vskip 10 pt
                \advance\subsecnum by 1 
                \defnum=1
                \leftline{\secfont \the\secnum.\the\subsecnum\ \quad #1}
                }

\def\definition{
                \advance\defnum by 1 
                \bf Definition 
\the\secnum .\the\defnum \rm \ 
                }

\def\lemma{
                \advance\lemmanum by 1 
                \par\bf Lemma  \the\secnum
.\the\lemmanum \rm \ 
                }

\def\theorem{
                \advance\lemmanum by 1 
                \par\bf Theorem  \the\secnum
.\the\lemmanum \rm \ 
               }

\def\cite#1{
				\secfont [#1]
				\rm
}

\vglue 20 pt

\centerline{\secfont QUATERNIONIC INVARIANTS of}
\centerline{\secfont VIRTUAL KNOTS and LINKS}

\bigskip

\bigskip

\centerline{\nam ANDREW BARTHOLOMEW, ROGER FENN${}^1$}
\centerline{\aff ${}^1$School of Mathematical Sciences, University of Sussex}
\centerline{\aff Falmer, Brighton, BN1 9RH, England}
\centerline{\aff e-mail addresses: rogerf@sussex.ac.uk, Andrewb@layer8.co.uk}
\centerline{\nam AMS classification 57M27}

\bigskip


\bigskip

\baselineskip=10 pt
\parskip=0 pt
\centerline{\nam ABSTRACT}
\leftskip=0.25 in
\rightskip=0.25in

{\nam In this paper we define and give examples of a family of polynomial
invariants of virtual knots and links. They arise by considering certain
2$\times$2 matrices with entries in a possibly non-commutative ring, for
example the quaternions. These polynomials are sufficiently powerful to 
distinguish the Kishino knot from any classical knot, including the unknot.

}

\leftskip=0 in
\rightskip=0 in
\baselineskip=13 pt
\parskip=\baselineskip

\noindent The contents of the paper are as follows
\smallskip

{\parindent=20pt

\parskip=0 pt
\item{1.}Introduction
\item{2.}Virtual Links
\item{3.}Switches: Definition and Examples
\item{4.}Labelling Diagrams
\item{5.}The Invariant $R$-module
\item{6.}The Invariant Polynomials
\item{7.}The Classical Case
\item{8.}Some Classifications of Switches
\item{9.}Some Calculations of the Polynomials
\item{10.}Questions and further developements
\item{11.}References
\item{12.}Tables

}
\parskip=\baselineskip

\section{Introduction}

This paper is organised as follows. In the next section we describe virtual
links and give reasons for studying them. In section 3 the conditions are
given for a 2$\times$2 matrix, $S$, to be a {\it linear switch} or {\it
switch} for short, see \cite{FJK}. The entries in the switch lie in some
associative ring, $R$, with identity which need not have commutative
multiplication. This allows us to define in section 5 an $R$-module for any
virtual link. The definition of the $R$-module comes from a labelling of a
diagram of the virtual link considered in section 4. If the ring allows
determinants then a sequence of polynomials in one or more variables can be
defined and this is treated in section 6. Virtual links also include classical
links but in that case these polynomials are constant as we see in section 7.
We would like to thank Hugh Morton for suggestions in this section and for
reading through an earlier draft. Finally various classification schemes of
switches are given together with tables. \np
\section{Virtual Links}

In this section we consider virtual links. For more details see \cite{K, FRS}
A diagram of a classical knot or link can be described by the Gauss code.
However not all Gauss codes can be realised as {\it classical} diagrams of
knots or links. Their realization may be dependant on the introduction of {\it
virtual crossings}. These are crossing which are neither above or below in
space but just indicate that the journey of the arc intersects the journey of
another arc. Virtual links are represented by oriented diagrams with ordinary
crossings as for classical knots and links together with these virtual
crossings. In addition to their application as a geometric realization of the
combinatorics of a Gauss code, virtual links have physical, topological and
homological applications. In particular, virtual links may be taken to
represent a particle in space and time which dissappears and reappears. A
virtual link may be represented, up to stabilisation, by a link diagram on a
surface. Finally an element of the second homology of a rack space can be
represented by a labelled virtual link, see \cite{FRS}. Since the rack spaces
form classifying spaces for classical links the study of virtual links may
give information about classical knots and links.

A {\bf diagram} for a virtual link is a 4-regular plane graph with extra
structure at its nodes representing the three types of crossings in the link.
A classical crossing of either sign is represented in the diagram in the usual
way. A virtual crossing is represented by two crossing arcs with a small
circle placed around the crossing point. The graph also lies implicitly on a
two-dimensional sphere $S^{2}$. {\bf Semi-arcs} go from one classical crossing
of the graph to another ignoring virtual crossings. This is distinct from a
classical link diagram where the {\it arcs} go from one undercrossing to
another.

Two such diagrams are {\bf equivalent} if there is a sequence of moves of the
types indicated in the figures below taking one diagram to the other. They are
the generalised Reidemeister moves and are local in character.

We show the classical Reidemeister moves as part (A) of Figure 1. These
classical moves are part of virtual equivalence where no changes are made to
the virtual crossings. Taken by themselves, the virtual crossings behave as
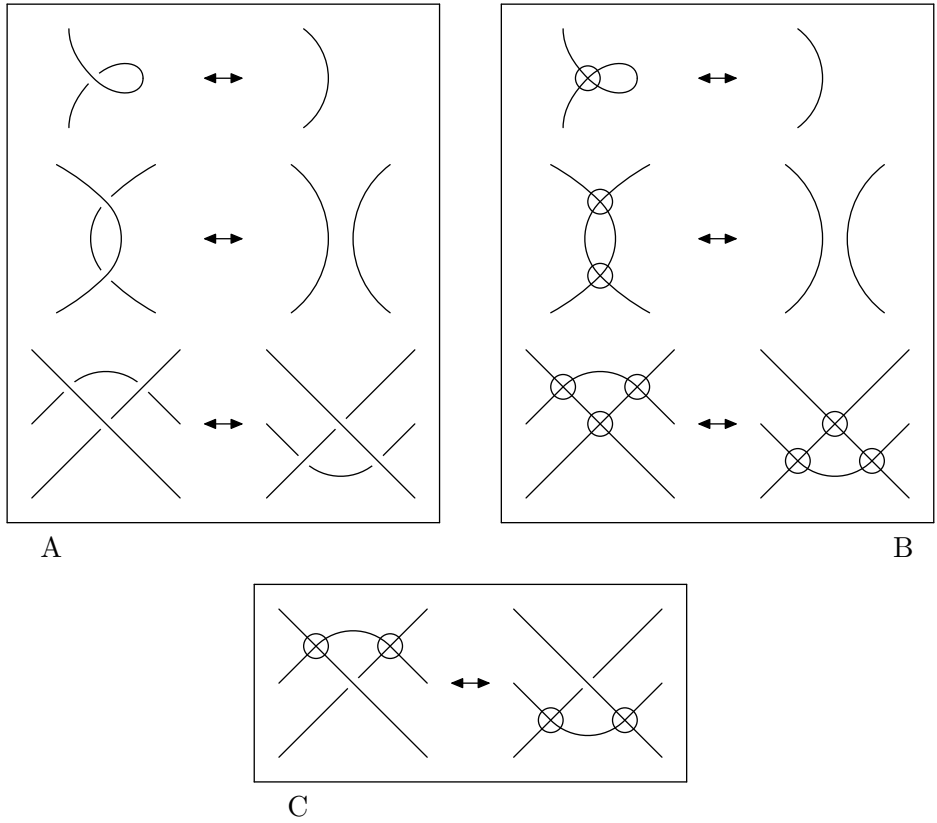
diagrammatic permutations. Specifically, we have the flat Reidemeister moves
(B) for virtual crossings as shown in Figure 1. In Figure 1 we also illustrate
a basic move (C) that interrelates real and virtual crossings. In this move an
arc going through a consecutive sequence of two virtual crossings can be moved
across a single real crossing. In fact, it is consequence of moves (B) and (C)
for virtual crossings that an arc going through any consecutive sequence of
virtual crossings can be moved anywhere in the diagram keeping the endpoints
fixed and writing the places where the moved arc now crosses the diagram as
new virtual crossings. This is shown schematically in Figure 2. We call the
move in Figure 2 the {\bf detour}, and note that the detour move is equivalent
to having all the moves of type (B) and (C) of Figure 1. This extended move
set (Reidemeister moves plus the detour move or the equivalent moves (B) and
(C)) constitutes the move set for virtual knots and links. \diagram \diagram

\section{Switches: Definition and Examples}
 
In this section we define the conditions needed for a 2$\times$2 matrix, $S$,
to be a {\it linear switch} or {\it switch} for short. The conditions, divided
into subsections, are invertability, Yang-Baxter and the existence of sideways
matrices. The reasons for these conditions should become clear in the next
section. Note that the definition of a switch in \cite{FJK} is more general.
\np
\subsection{Inverting a 2$\times$2 Matrix}

Suppose $S$ is the $2\times 2$ matrix with entries in a ring $R$,
$$S=\pmatrix{A & B \cr C & D\cr}.$$ The ring $R$ is associative and has a
multiplicative identity element 1 but need not have commutative
multiplication. We call an element of the ring a {\bf unit} if it has a two
sided multiplicative inverse.

The proof of the following lemma may be safely left with the reader.

\lemma{The matrix $S$ is a unit in the ring of 2$\times$2 matrices with
entries in $R$ if either
$B, D$ and $\Delta=B^{-1}A-D^{-1}C$ are units or $A, C$ and 
$\Delta'=C^{-1}D-A^{-1}B$ are units. In the first case
$$S^{-1}=\pmatrix{\Delta^{-1}B^{-1} & -\Delta^{-1}D^{-1} \cr
-D^{-1}C\Delta^{-1}B^{-1} & B^{-1}A\Delta^{-1}D^{-1}\cr}.$$
In the second case
$$S^{-1}=\pmatrix{C^{-1}D\Delta'^{-1}A^{-1} & -A^{-1}B\Delta'^{-1}C^{-1} \cr
-\Delta'^{-1}A^{-1} & \Delta'^{-1}C^{-1}\cr}.$$
}
\qed

The hypothesis for the above lemma is sufficient but not quite necessary. For
example $\pmatrix{1&\alpha\cr 0&1}$ is invertible even if $\alpha$ is not.
Note that the condition for an inverse to exist can also be written
$B^{-1}A\ne D^{-1}C$ if the ring is a division ring.

The theory of inverting matrices with entries in a non-commutative ring can
differ from the commuting case. As an example consider the following $2\times
2$ matrix with entries in the quaternions $\H$. Let $$S=\pmatrix{1& k\cr i &
j\cr}\hbox{ then }S^{-1}=\pmatrix{1/2& -i/2\cr -k/2 & -j/2\cr}.$$ The
``obvious'' determinant $j-ki$ is zero. Moreover the transpose of $S$ has no
inverse since then $\Delta=0$.

\subsection{The Yang-Baxter Equations}

The next condition to consider is
$$ (S\times id)(id\times
S)(S\times id)=(id\times S)(S\times id)(id\times S).$$
Here $S\times id$ and $id\times S$ are the 3$\times$3 matrices
$$S\times id=\pmatrix{A&B&0\cr C&D&0\cr 0&0&1\cr}\quad
id\times S=\pmatrix{1&0&0\cr 0&A&B\cr 0&C&D\cr}.$$
These equations are a specialization of the set theoretic Yang-Baxter
equations considered by  V. Drinfeld in \cite{Dr}.

The Yang-Baxter equations imply
the seven equations $$\matrix{%
1: A=A^2+BAC\hfill  &2: [B,A]= BAD\hfill \cr 
3: [C,D]= CDA\hfill  &4: D=D^2+CDB\hfill  \cr 
5: [A,C]= DAC\hfill  &6: [D,B]= ADB \hfill \cr 
\hfill 7: [C,B]=& ADA-DAD \hfill \cr}$$
where $[X,Y]$ denotes the commutator $XY-YX$.

\lemma{Assume that either $A, B, C$ or $B, C, D$ are units and the seven
equations above are satisfied. Then $A-1, D-1$ are also units and the seven
equations can be reduced to the first four equations. By a further refinement
these equations can be reduced to two equations in two unknowns.}

{\bf Proof } Assume that $A, B, C$ are units. Firstly we write $C, D$ in terms
of $A, B$ using equations 1 and 2. So $$\matrix{%
C=&A^{-1}B^{-1}A-A^{-1}B^{-1}A^2=A^{-1}B^{-1}A(1-A) &
D=&1-A^{-1}B^{-1}AB\cr}.$$ We see that 5 is easily satisfied if we substitute
these values and clearly $A-1$ and $D-1$ are units.

Now look at equation 4. The right hand side minus the left hand side is
$A^{-1}B^{-1}A\Theta B$ where
$$\Theta=BA^{-1}B^{-1}A-A^{-1}B^{-1}AB-A+B^{-1}AB.$$

The same difference for equation 6 is $-\Theta B$ and for equation 7 is
$\Theta(1-A)$. So $\Theta=0$ implies equations 6 and 7. The converse is also
true since $B, A-1$ are units.

It is also easily seen that $(D-1)B^{-1}(A-1)C^{-1}=1$ and we note that for
future use.

Since $C,\ D$ can be eliminated using equations 1 and 2 we
are finally left with two equations, in two unknowns $A$ and $B$. 

A symmetric argument works if $B, C, D$ are units.\qed

The number of equations can be reduced to one, see \cite{B,F}.

\subsection{The Final Definition and Sideways Matrices}

Summing up, a matrix $S=\pmatrix{A & B \cr C & D\cr}$ is a {\bf linear
switch} if 

{\parindent=20pt
1. $B, C$ are units and either $D$ and $\Delta=B^{-1}A-D^{-1}C$ are units or
$A$ and $\Delta'=C^{-1}D-A^{-1}B$ are units.

2. The four Yang-Baxter equations
$$\matrix{%
1: A=A^2+BAC\hfill  &2: [B,A]= BAD\hfill \cr 
3: [C,D]= CDA\hfill  &4: D=D^2+CDB\hfill  \cr}
$$ are satisfied.
}

The {\bf sideways matrices} $S^+_-$ and $S^-_+$ are defined  by
$$S^+_-=\pmatrix{%
DB^{-1}&C-DB^{-1}A\cr B^{-1}&-B^{-1}A\cr}\quad
S^-_+=\pmatrix{%
-C^{-1}D&C^{-1}\cr B-AC^{-1}D&AC^{-1}\cr}
$$
Note for future reference that both $S^+_-$ and $S^-_+$ are invertible if $S$
is and that $(S^{-1})^+_-=(S^-_+)^{-1}$ and $(S^{-1})^-_+=(S^+_-)^{-1}$.
Also 
$$S^+_-(a,a)=(\lambda a,\lambda a)\hbox{ and }
S^-_+(a,a)=(\lambda^{-1} a,\lambda^{-1} a)$$
where $\lambda=B^{-1}(1-A)=(1-D)^{-1}C$. So the sideways matrices preserve
the diagonal. This has the curious consequence that a linear switch which is
a birack is also a biquandle in the sense of \cite{FJK}.

It is an easy exercise to verify that the only switches with a commutative
ring are $$\pmatrix{ 1-\mu\lambda& \mu \cr\lambda & 0\cr} \quad\pmatrix{0 &
\mu \cr \lambda & 1-\mu\lambda\cr}$$ for some units $\lambda,\ \mu$ in the
ring.

Either is called the {\bf Alexander} switch, see \cite{FJK}.
Note that
$$\pmatrix{1&0\cr 0&\lambda^{-1}\cr}\pmatrix{1-\mu\lambda&\mu\cr\lambda&0\cr}
\pmatrix{1&0\cr 0&\lambda\cr}=\pmatrix{1-\mu\lambda&\mu\lambda\cr 1&0\cr}.$$
The matrix on the right is a version of the Burau matrix.

The identity is invertible and a solution of the Yang-Baxter equations but
is not a linear switch because $B, C$ are zero.

For a non-commutative example assume that the ring is the quaternion division
algebra $\H$. Then
$$S=\pmatrix{%
1+i&-j\cr j&1+i\cr}$$
is called the {\bf Budapest switch}. It is only one of many such 
solutions found by analysis and a computer search.

It is easily checked that if $$S=\pmatrix{A & B \cr C & D\cr}$$ is a linear
switch then so is $$S(t)=\pmatrix{A & tB \cr t^{-1}C & D\cr}$$ where $t$ is
any variable in the centre of the ring. By analogy with the quaternionic case
we call $t$ a {\bf real} variable. Further switches are given by $S^{-1}$,
because $S\to S^{-1}$ is an antiautomorphism and $S^\dagger=\pmatrix{D&C\cr
B&A\cr}$ is also a switch, by inspection. If the entries of $S$, a linear
switch, are quaternions, then $S^*=\pmatrix{\overline{A}&\overline{C}\cr
\overline{B}& \overline{D}\cr}$ is a solution, because hermitian involution is
also an antiautomorphism, and so is
$S^{\dagger*}=\pmatrix{\overline{D}&\overline{B}\cr \overline{C}&
\overline{A}\cr}$.

\section{Labelling Diagrams}

In this section, given a switch $S$ with entries in $R$, we define a {\it
labelling} or colouring, $\cal{L}$, of the semi-arcs of a virtual link
diagram, $D$, by elements of $R$ in such a way that after a Reidemeister move
converting $D$ into $D'$ there is a uniquely defined labelling ${\cal L}'$ of
$D'$ which is unchanged outside of the disturbence caused by the Reidemeister
move. It follows that if $D_1$ and $D_2$ are diagrams representing the same
virtual link and $D_1\to\cdots\to D_2$ is a sequence of Reidemeister moves
transforming $D_1$ into $D_2$. Then any labelling ${\cal L}_1$ of $D_1$ is
transfered via the sequence of Reidemeister moves to a labelling ${\cal L}_2$
of $D_2$. In particular the set of labellings of $D_1$ is in bijective
correspondence with the set of labellings of $D_2$, albeit not by a uniquely
defined bijection.

Let the edges of a positive real crossing in a diagram be arranged diagonally
and called geographically {\it NW, SW, NE} and {\it SE}. Assume that initially
the crossing is oriented and the edges oriented towards the crossing from left
to right ie west to east. The {\bf input} edges, oriented towards the
crossing, are in the west and the edges oriented away from the crossing, the
{\bf output} edges, are in the east. Let $R$ be a labelling set and let $a$
and $b$ be labellings from $R$ of the input edges with $a$ labelling SW and
$b$ labelling NW. For a positive crossing, $a$ will be the label of the
undercrossing input and $b$ the label of the overcrossing input. Suppose now
that $S(a,b)^T=(c,d)^T$. Then we label
the undercrossing output NE by $d$ and we label the overcrossing output SE by
$c$. 

For a negative crossing the direction of labelling is reversed. So $a$ labels
SE, $b$ labels NE, $c$ labels SW and $d$ labels NW.

Finally for a virtual crossing the labellings carry across the strings. This
corresponds to the twist function $T(a,b)=(b,a)$. Manturov in \cite{M} has
generalized this to $T(a,b)=(\epsilon b,\epsilon^{-1}a)$. However as we shall
see later this generalization leads to the same polynomial considered by
Silver and Williams, \cite{SW}.
\np

The following figure shows the labelling for the three kind of
crossings.
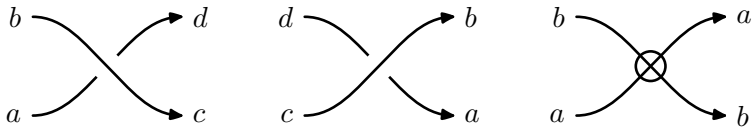
\diagram
\centerline{$c=Aa+Bb\quad d=Ca+Db$}

It is convenient to think of the action from left to right on a positive
crossing as being the action of $S$, the action from right to left as being
$S^{-1}$, the action from top to bottom as being $S^-_+$ and the action from
bottom to top as being $S^+_-$. For a negative crossing the actions are equal
but with opposite orientation.
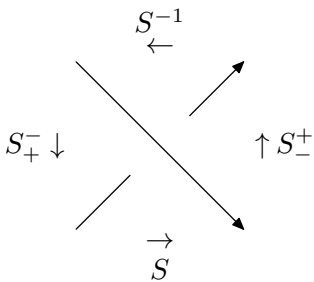
\diagram

\theorem{The set of labellings of two diagrams representing the same virtual
link are in bijective equivalence.}

{\bf Proof } Full details can be found in \cite{FJK}. It should now not be
difficult for the reader to see how a labelling on a diagram is extended to a
labelling on the result of a Reidemeister move. For example consider a
Reidemeister 1 move which introduces a new classical crossing. Then the new
semi-arc needs a consistent label. But this happens because the sideways
matrices preserve the diagonal. A Reidemeister 2 move either extends the
labelling because $S$ is invertible or because one of the sideways operators
are invertible. A Reidemeister 3 move labelling follows from the Yang Baxter
equations. \qed

\section{The Invariant $R$-module}

Assume for simplicity that we are dealing with a knot. The link case is
similar and details can safely be left to the reader.

The set of labellings of a diagram by the ring $R$ with fixed switch $S$ is a
right $R$-module and we can give a finite (square) presentation.

More precisely let $D$ be a diagram of a virtual link with $n$ classical
crossings. The {\bf semi-arcs} of the diagram run from one classical crossing
to the next. The virtual crossings are ignored. Following the orientation of
the knot, label the semi-arcs with $R$-variables $x_1, x_2, \ldots, x_{2n}$.
By an $R$-variable we mean a symbol standing in for any element of $R$.

At each crossing there is a relation of the form
$$\pmatrix{A&B\cr C&D\cr}\pmatrix{x_i\cr x_j\cr}=
\pmatrix{x_{j+1}\cr x_{i+1}\cr}
\hbox{ or }
\pmatrix{A&B\cr C&D\cr}\pmatrix{x_i\cr x_j\cr}=
\pmatrix{x_{j-1}\cr x_{i-1}\cr}$$
depending on whether the crossing is positive or negative. As is the usual
custom, indices are taken modulo $2n$.

The relations can now be written in matrix form as $M{\bf x}={\bf 0}$ where
$M$ is a $2n\times 2n$ matrix and ${\bf x}=(x_1,x_2 \ldots, x_{2n})^T$. The
non-zero entries in each row of the matrix are $A, B, -1$ or $C, D, -1$.

Let ${\cal M}={\cal M}(S,D)$ be the module defined by these relations. We now
show that the modules defined by diagrams representing the same virtual link
are isomorphic. We do this by showing that a single Reidemeister move defines
an isomorphism. The proof has the same structure as the proof, say, that the
Alexander module of a classical link is an invariant as in \cite{A} but we
give the details because of the care needed due to non-commutativity.

\theorem{The module ${\cal M}$ defined above is invariant under the
Reidemeister moves and is therefore a virtual link invariant.}

{\bf Proof} 
Any module defined by a presentation of the form $M{\bf x}={\bf 0}$ is
invariant under the following moves and their inverses applied to the matrix
$M$. 
 
{\parindent=20pt
\item{1.} permutations of rows and columns,
\item{2.} multiplying any row on the left or any column on the right by a
unit,
\item{3.} adding a left
multiple of a row to another row or a right multiple of a column to another
column, 
\item{4.} changing $M$ to $\pmatrix{x&{\bf u}\cr{\bf 0}&M\cr}$ where
$x$ is a unit, ${\bf u}$ is any row vector and ${\bf 0}$ is a zero column
vector,
\item{5.} repeating a row.
}

An {\bf elementary} matrix of type 1 is a permutation matrix. An elementary
matrix of type 2 is the identity matrix with one diagonal entry replaced by a
unit and an elementary matrix of type 3 is a square matrix with zero entries
except for $1$'s down the diagonal and one other entry off diagonal.

The operations $i.$ above for $i=1, 2, 3$ are equivalent to multiplying $M$ on
the right or left by an elementary matrix of type $i$.

For example any switch can be written as a product of elementary matrices as
follows
$$\pmatrix{%
A & B \cr C & D\cr}=
\pmatrix{%
A & 0 \cr 0 & 1}
\pmatrix{%
1 & 0 \cr C & 1\cr}
\pmatrix{%
1 & 0 \cr 0 & C\Delta'\cr}
\pmatrix{%
1 & A^{-1}B \cr 0 & 1}
$$
if $A$ and $C\Delta'$ are units and a similar formula if $D$ and $B\Delta$
are units.

Now consider the module ${\cal M}$ defined above. Clearly the presentation
is unaltered by any of the basic moves which involve the virtual crossing. So
we look to see the changes induced by the classical Reidemeister moves and
check that the presentation matrix is only changed by the above 5 moves.

Firstly, consider a Reidemeister move of the first kind. 
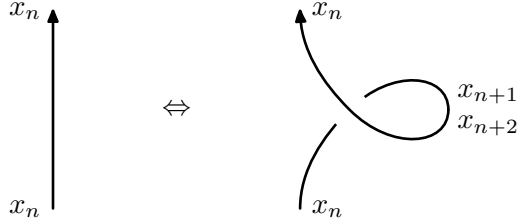
\diagram This
introduces (or deletes) two new equal generators $x_{n+1}=x_{n+2}$. Because
$S^-_+$ and $S^+_-$ preserve the diagonal, (the biquandle condition, see
\cite{FJK}) the output ($x_{n+3}$) is the same as the input ($x_{n}$). The
generator $x_{n+1}$ is equal to $\lambda^{-1}x_{n}$ where
$\lambda=B^{-1}(1-A)$.

So up to reordering of the columns the relation matrix is changed by
$$M\ \Leftrightarrow\ \pmatrix{M&\matrix{{\bf 0}&{\bf 0}\cr}\cr
\matrix{{\bf 0}&0\cr {\bf 0}&1\cr}&\matrix{1&-1\cr\lambda&0\cr}\cr}
\sim\pmatrix{M&\matrix{{\bf 0}&{\bf 0}\cr}\cr
\matrix{{\bf 0}&0\cr {\bf 0}&1\cr}&\matrix{0&-1\cr\lambda&0\cr}\cr}$$
Since $\lambda$ is a unit this does not alter the module.

There are other possible inversions and mirror images of the above which
can be dealt with in a similar fashion.

Secondly, consider a Reidemeister move of the second kind. 
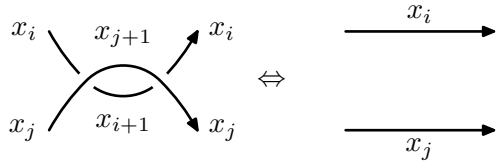
\diagram Again the
outputs are unchanged from the inputs $x_j, x_i$ because of the relation
$S^{-1}S=1$.

Two new generators $x_{i+1}$ and $x_{j+1}$ are introduced (or deleted).
They are related by the equations
$$x_{j}=Ax_{i+1}+Bx_{j+1}\hbox{ and }x_{i}=Cx_{i+1}+Dx_{j+1}.$$

This has the following effect on the relation matrix.
$$M\ \Leftrightarrow\ \pmatrix{M&\matrix{{\bf 0}&{\bf 0}\cr}\cr
\matrix{{\bf 0}&0&-1\cr{\bf 0}&-1&0\cr}&\matrix{A& B \cr C& D\cr}\cr}.
$$
Since $S$ is a product of elementary matrices this does not alter the module.

The other possible inversions and mirror images of the above can be dealt with
in a similar fashion but it is worth looking at the case where the two arcs
run in opposite directions. The right outputs are unchanged from the left
inputs by the relation $S^+_-(S^+_-)^{-1}=1$.
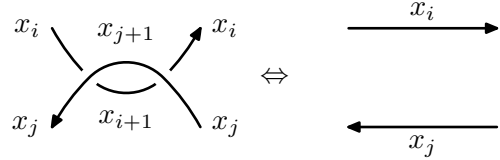
\diagram
The changes to the
relation matrix are given by $$M\ \Leftrightarrow\ \pmatrix{M&\matrix{{\bf
0}&{\bf 0}\cr}\cr \matrix{{\bf 0}&-1&A\cr{\bf 0}&0&C\cr}&\matrix{B& 0
\cr D&-1\cr}\cr}\sim
\pmatrix{M&\matrix{{\bf
0}&{\bf 0}\cr}\cr \matrix{{\bf 0}&-1&A\cr{\bf 0}&0&0\cr}&\matrix{B& 0
\cr 0&-1\cr}\cr} $$ 

The module is unaltered because $B$ is a  unit.

Finally, consider a Reidemeister move of the third kind. 
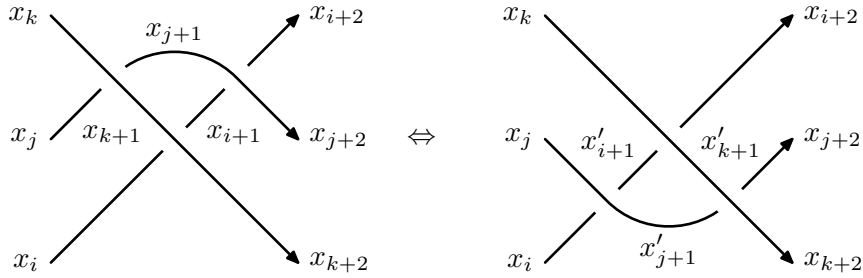
\diagram The outputs
$x_{i+2}, x_{j+2}, x_{k+2}$ are unaltered by the Reidemeister move because of
the Yang-Baxter equations. The inner generators $x_{i+1}, x_{j+1}, x_{k+1}$
are related to the inputs $x_{i}, x_{j}, x_{k}$ by the following matrix
$$\pmatrix{C&DA&DB\cr 0&C&D\cr 0&A&B\cr}$$ and the inner
generators $x'_{i+1}, x'_{j+1}, x'_{k+1}$ are related to $x_{i}, x_{j}, x_{k}$
by the following matrix $$\pmatrix{C&D&0\cr A&B&0\cr
AC&AD&B\cr}.$$ Both are the product of elementary matrices and the
proof is finished. \qed

\section{The Invariant Polynomials}

We have seen that associated to every linear switch $S$ and virtual link
diagram $D$ there is an associated $R$-module, ${\cal M}$, which is a link
invariant. Moreover this module has a finite presentation of the form $M{\bf
x}={\bf 0}$. If the ring is commutative then there is a sequence of invariants
of ${\cal M}$, the {\bf elementary ideals}, $E_i$ which can be defined from
the presentation. The details can be found in the book by Crowell and Fox
\cite{CF}. Suppose an $R$-module is defined by a matrix $M$, (not necessarily
square but with $c$ columns and $r$ rows, $r\ge c$). Let $E_0$ be the ideal
generated by the largest possible sub-determinants of $M$, that is those of
size $c\times c$. Let $E_1$ be the ideal generated by the $c-1\times c-1$
sub-determinants of $M$ and so on. The elementary ideals form an
ascending chain of ideals $$\{0\}\subset E_0\subset E_1\subset\cdots\subset
E_c=E_{c+1}=\cdots=R.$$ The ideal $E_{c-1}$ is generated by the elements of
$M$.

In our situation the matrix $M$ is square so $E_0$ is principle, generated by
the determinant of $M$. This has the following interpretation. Let $\adj M$
denote the adjugate matrix of $M$ consisting of the codimension 1 minors of
$M$. Let $(x_1,\ldots,x_n)$ be the generators of $\cal M$. Multiplying the
defining equations on the left by $\adj M$ yields the equation
$$\det(M)(x_1,\ldots,x_n)=0.$$ It follows that if $(x_1,\ldots,x_n)\ne0$ then
either $\det(M)=0$ in the ring or ${\cal M}$ has torsion.

Assume now that $R$ is a gcd-ring. That is $R$ is an integral domain such that
any set of elements has a gcd which is well defined up to multiplication by a
unit. Examples of gcd-rings are the (Laurent) polynomial rings
$p[t_1^{\pm1},\ldots,t_k^{\pm1}]$ where $p$ is a field or the integers $\Z$.
It follows that we can define the $i$-th {\bf ideal polynomial} by the rule
$$\Delta_i({\cal M})=\hbox{gcd}\{E_i\}.$$ Then $\Delta_i({\cal M})$ is the
generator of the smallest principle ideal containing $E_i$. We have a chain of
divisions $$1=\cdots=\Delta_n|\Delta_{n-1}|\cdots|\Delta_0=\det(M)|0.$$

Let us illustrate this with the Alexander switch. Here
$R=\Z[\lambda,\lambda^{-1},\mu,\mu^{-1}]$ and the switch is defined by the
matrix $$\pmatrix{0 & \mu \cr \lambda & 1-\mu\lambda\cr}.$$ We shall call
$\det(M)=\Delta_0(K)$ the 0-th {\bf Alexander polynomial} where $K$ is the
class of the diagram $D$. For classical knots $K$, it is always the case that
$\Delta_0(K)=0$. However this need not be true for virtual knots or links. For
example consider the {\bf virtual trefoil} as shown in the figure. If we label
as indicated then the module has a presentation with 3
generators $a, b, c$ and relations $c=\mu a,\ a=\mu\lambda b+
\mu(1-\mu\lambda)a,\
b=\lambda c+(1-\mu\lambda)(\lambda b+(1-\mu\lambda)a)$.
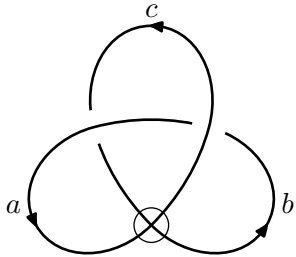
\diagram
If we eliminate $c=\mu a$ we arrive at the following equations.

$$\matrix{(\mu-\lambda\mu^2-1)a&\hfill+\lambda\mu
b=&0\cr (\lambda^2\mu^2-\lambda\mu+1)a&\hfill+(\lambda-\lambda^2\mu-1)b=&0\cr}$$

The determinant of these equations is
$$\Delta=(\lambda-1)(\lambda\mu-1)(\mu-1).$$
Note that the fundamental rack
(and hence group defined by the Wirtinger relations) is trivial.

Consider the Kishino knots $K_1, K_2$ and $K_3$ illustrated below.
$$\ddiagram\qquad\ddiagram\qquad\ddiagram$$
All are ways of forming the connected sum of two
unknots. Under the symmetry which changes positive crossings into negative
crossings and vice-versa leaving virtual crossings alone, $K_1$ is transformed
into $K_2$ and $K_3$ is invariant. Both have trivial racks and Jones
polynomial. The Alexander polynomial $\Delta_0$ is zero in all three cases. On
the other hand for $K_1$, the 1st Alexander polynomial, $\Delta_1$ is
$1+\mu-\lambda\mu$ and for $K_2$, $\Delta_1$ is $1+\lambda-\lambda\mu$. Since
these are neither units nor associates in the ring, $K_1, K_2$ are non trivial
and non amphich\ae ral in the above sense.  Since $\Delta_1$ in the case of $K_1, K_2$ is not an
Alexander polynomial in the form $f(\lambda\mu)$ it follows that $K_1, K_2$
are not equivalent to classical knots, see \cite{S,W}.

The 1st Alexander polynomial $\Delta_1$ of $K_3$ is 1. So we need new
invariants to distinguish $K_3$ from the trivial knot. The Alexander switch is
the only commutative case so we need a non-commutative ring. In the case of
the entries in the switch this means that the pairs $(A, B)$, $(A, C)$, $(D,
B)$, $(D, C)$, may not commute although other pairs may.

Firstly we need a definition of the determinant of a matrix with entries in a
non-commutative ring such as the quaternions. There are various definitions of
determinants in this case and we refer to \cite{As} for details. For our
purposes let $M_n(R)$ denote the ring of $n\times n$ matrices with entries in
the ring $R$. We will only consider $R=\R,\ \C,\ \H$, the reals, complex
numbers and quaternions respectively or polynomials with coefficients in these
rings.

Any quaternion $Q$ can be written as a pair of complex numbers $a, b$ by the
formula $Q=a+bj$. There is an injective ring homomorphism $\psi: M_n(\H)\to
M_{2n}(\C)$ given by $$\psi(a+bj)=\pmatrix{a&b\cr
-\overline{b}&\overline{a}\cr}.$$ Define the {\bf Study} determinant by
$$\d(M)=\det(\psi(M)).$$ Note that the definition is slightly different from
the one given in \cite{As} but its properties are the same.

This Study determinant has the following properties
{\parindent=20pt
\item{1.} $\d(M)=0$ if and only if $M$ is singular.
\item{2.} $\d(MN)=\d(M)\d(N)$
\item{3.} $\d(M)$ is unaltered by adding a left multiple of a row to another
row or a right multiple of a column to another column.
\item{4.} $\d\pmatrix{x&{\bf u}\cr{\bf 0}&M\cr} =|x|^2\d(M)$ where ${\bf u}$
is any row vector and ${\bf 0}$ is a zero column vector.
\item{4'.} $\d\pmatrix{x&{\bf 0}\cr{\bf v}&M\cr} =|x|^2\d(M)$ where ${\bf v}$
is any column vector and ${\bf 0}$ is a zero row vector.
\item{5.} $\d(M^*)=\d(M)$ where $M^*$ denotes Hermitian conjugate.
\item{6.} If $M$ is non-singular then $\d(M)$ is a positive real number.
\item{7.} $\d(M)$ is unaltered by permutations of rows and columns.
}

We now define the elementary ideals in the non-commutative case. As before,
let $E_0$ be the ideal generated by the largest possible sub-determinants of
$M$. Let $E_1$ be the ideal generated by the second largest possible
sub-determinants of $M$ {\it plus} $E_0$ and so on. This extra condition is
needed in order that the elementary ideals form an ascending chain and we
cannot expand by rows or columns as in the case of commutative rings.

It is not hard to see that these ideals are invariant under the 5 moves above.

Let $S$ be a switch with quaternion entries. We can make the entries lie in
$\H[t,t^{-1}]$ by replacing $S$ with $S(t)$. The polynomials
$\Delta_{i}^{\H}(D,S)$ are defined as in the commutative case. That is
$$\Delta_{i}^{\H}(D,S)=\hbox{gcd}\{E_i\}.$$ By the properties of the Study
determinant they will be real (Laurant) polynomials with roots in
complex conjugate pairs and are defined up to multiplication by a unit.

We can normalise the polynomials $\Delta_{i}^{\H}(D,S)$ by multiplying by a
suitable unit $\pm t^n$ so that we get a genuine polynomial with a positive
constant term.


For a worked calculation of $\Delta_{0}^{\H}$, consider again the diagram of
the virtual trefoil above.

The equations on the generators $a, b, c$ are
$$\eqalign{
tBa+Ab-c=&0\cr ADa+(t^{-1}AC-1)b+tBc=&0\cr (t^{-1}CD-1)a+t^{-1}C^2b+DAb=&0\cr}
.$$
Eliminating $c=tBa+Ab$ the relation matrix is
$$\pmatrix{
AD+t^2B^2&t^{-1}AC+tBA-1\cr t^{-1}CD+tDB-1& t^{-1}C^2+DA\cr}.$$
If we now take the values  $A=1+i, B=-j, C=j, D=1+i$ of the Budapest switch,
the normalised determinant is $\Delta_0^{\H}=1+2t^2+t^4$.

For the three Kishino knots $\Delta_0^{\H}=0$ and
$\Delta_1^{\H}=1+(5/2)t^2+t^4$ showing in particular that $K_3$ is
non-trivial.

The polynomials $\Delta_i^{\H}$ are defined by deleting rows and columns of a
presentation matrix $M$ and taking the determinant of the image under the
automorphism $\psi$ which now has complex entries. But we can reverse this
procedure and obtain a new set of polynomials $\Delta_i^{\C}$. The even order
polynomials $\Delta_{2i}^{\C}$ will divide $\Delta_i^{\H}$.

We can continue to the real case under the well known embedding of the complex
numbers in the ring of $2\times 2$ matrices. The resulting polynomials are
denoted by $\Delta_i^{\R}$.

At the moment we have not calculated any examples but they may provide useful
invariants.

\section{The Classical Case}

This section suggests that these polynomials do not provide more information
about classical knots and links. However it is a method of distinguishing
virtual knots from classical ones. The idea behind the proofs is to use the
fact that any classical link is the closure of a braid, \cite{B}. This result
can easily be extended to virtual links.

{\bf Virtual} braids are defined in a similar manner to classical braids. We
consider a diagramatic representation of a virtual braid as a set of $n$
strings travelling from a vertical line of $n$ points on the left to a
translated copy of these points on the right. The strings are to be monotonic
in the left to right motion. They may cross as in the classical case with a
right handed cross (corresponding to $\sigma_i$) or they may incur a left
handed cross (corresponding to $\sigma_i^{-1}$) or they may have a virtual 
crossing. The virtual crossing is indicated algebraically by $\tau_i$.

The group of virtual braids on
$n$ strings is denoted by $VB_n$. Kamada's presentation of the virtual braid
group is as follows, \cite{KK}.

Let generators of $VB_n$ be $\sigma_1, \ \sigma_2, \ ..., \sigma_{n-1}$ where
$\sigma_i$ corresponds to the real positive crossing of the $i$-th and $(i+1)$-th
string and $\tau_1, \ \tau_2, \ ..., \ \tau_{n-1}$ where $\tau_i$ corresponds
to the virtual crossing of the $i$-th and $(i+1)$-th string. The following
relations hold:

{\parindent=20pt

\item{i)}Braid relations $$\matrix{ \sigma_i \sigma_j= \sigma_j \sigma_i
,\qquad |i-j|>1 \cr \sigma_i \sigma_{i+1} \sigma_i=\sigma_{i+1} \sigma_i
\sigma_{i+1}\hfill \cr } $$ \item{ii)}Permutation group relations $$\matrix{
{\tau_i}^2=1 \hfill \cr \tau_i \tau_j= \tau_j \tau_i ,\qquad |i-j|>1 \cr
\tau_i \tau_{i+1} \tau_i=\tau_{i+1} \tau_i \tau_{i+1}\hfill \cr } $$
\item{iii)}Mixed relations $$\matrix{ \sigma_i \tau_j= \tau_j \sigma_i ,\qquad
|i-j|>1 \cr \sigma_i \tau_{i+1} \tau_i=\tau_{i+1} \tau_i \sigma_{i+1}\hfill
\cr } $$}
At this stage the reader may care to compare these relations with the
relations for the braid-permutation group, see \cite{FRR}. The
braid-permutation group is a quotient of the virtual braid group. The
corresponding knot-like objects are called {\bf welded} knots.

Any linear switch $S$ with entries in the ring $R$ defines a representation of
the braid group $B_n$ into the group of invertible $n\times n$ matrices with
entries in $R$ by sending the standard generator $\sigma_i$ to
$S_i=(id)^{i-1}\times S \times (id)^{n-i-1}$. Denote the representation by
$\rho=\rho(S, n)$. This can be extended to virtual braids by sending $\tau_i$
to $T_i=(id)^{i-1}\times T \times (id)^{n-i-1}$ where $T=\pmatrix{0&1\cr
1&0\cr}$. Denote the extended representation by $\rho'=\rho'(S, T, n)$.

\theorem{1. For all linear switches $S$ the representation, $\rho(S, n)$, of
the braid group $B_n$ is equivalent to $\rho(S(t), n)$.

2. There exists a non-zero vector ${\bf z}=(z_0, z_1, \ldots, z_{n-1})^T$ such
that $S_i({\bf z})={\bf z}$. In particular the representation $\rho$ is
reducible.

3. The representation, $\rho'(S(t), T, n)$, of the virtual braid group,
$VB_n$, is equivalent to \hfil\break $\rho'(S, T(t), n)$
where  $T(t)=\pmatrix{0&t\cr t^{-1}&0\cr}$.}

{\bf Proof } Let $\Lambda=\diag\{1, t, t^2, \ldots, t^{n-1}\}$. Then
$\Lambda^{-1}S_i\Lambda=S(t)_i$ and $\Lambda^{-1}T_i\Lambda=T(t)_i$.

The vector ${\bf z}$ is defined by ${\bf z}=(1,\lambda, \lambda^2, \ldots,
\lambda^{n-1})^T$ where $\lambda=B^{-1}(1-A)$. \qed

This shows that there is no gain in the representation of $B_n$ by replacing
$S$ with $S(t)$. It also shows that Manturov's invariant in \cite{M} is
equivalent to the generalised Alexander polynomial.

\theorem{For all classical knots and all switches where the polynomials are
defined, $\Delta_i$ is constant (independant of $t$) and in particular
$\Delta_0=0$.}

{\bf Proof } Consider a diagram where the link is the closure of a classical
braid. Suppose that the switch is $S(t)$. Then the defining matrix is of the
form $M-I$ where $M$ is the representing matrix of the braid group and $I$ is
the identity matrix. By the above $\Lambda(M-I)\Lambda^{-1}$ has switch
$S(1)$.

Moreover $(M-I){\bf z}=0$ so any determinant of $M-I$ will be zero. \qed

As a consequence of this theorem it follows that the Kishino knot $K_3$ is not
a classical knot. It also verifies, see \cite{SW}, that the Silver-Williams
version of the Alexander polynomial is just the original for classical knots
in the form $\Delta(\lambda\mu)$.

We illustrate this result by calculating this constant for the Budapest
switch. Recall that the determinant of a knot is the value of the Alexander
polynomial evaluated at $-1$.

\theorem{For the Budapest switch,
$S=\pmatrix{ 1+i&-j\cr j&1+i\cr}$ and for any classical knot,
$\Delta_1$ is the square of the knot determinant.}

{\bf Proof } Define the $n\times n$ lower triangular matrix ${\cal B}_n$
inductively as follows. Put ${\cal B}_1=1$ and suppose that ${\cal B}_{n-1}$
is defined with bottom row $x_1,x_2,\ldots,x_{n-1}$. Then ${\cal B}_n$ is
obtained from ${\cal B}_{n-1}$ by adjoining the bottom row
$y_1,y_2,\ldots,y_{n}$ where $$y_i=kx_i,\ i=1,2,\ldots,n-2,\
y_{n-1}=(k-j)x_{n-1}\hbox{ and }y_n=jx_{n-1}.$$ So for example $${\cal
B}_3=\pmatrix{1&0&0\cr k-j&j&0\cr i-1& 1-i&-1\cr}.$$ Then ${\cal
B}_n^{-1}S_i{\cal B}_n=D_i$ where $D=\pmatrix{0&1\cr -1&2\cr}$, which is the
Burau matrix with variable $-1$. The Study determinant now gives
the square of the Alexander polynomial evaluated at $-1$. \qed

It follows from the proof above that the representation of the classical
braid group when $S$ is the Budapest switch is equivalent to a variant of
the Burau representation. We conjecture that this is always the case.
\section{Some Classifications of Switches}

Consider now $2\times2$ linear switches with quaternion entries. We have
already seen that a switch $S$ also defines the switch $S(t)$.

The next lemma looks at the constraints on the entries of $S$. We can write
any quaternion $Q\in\H$ as the sum of a real part ${\cal R}(Q)\in \R$, and a
purely quaternionic part ${\cal P}(Q)\in \R^3$. Let $S^3$ denote the group of
unit quaternions and let $S^2$ denote the set of pure unit quaternions. So
$S^2$ is $\sqrt{-1}$. Let $\C$ be the complex numbers in $\H$ and let $S^1$
denote the group of unit complex numbers.

\lemma{Let $S=\pmatrix{A & B \cr C & D\cr}$, with quaternion entries, be a
switch. Then, in the non-commutative case, $|A|<2$ and ${\cal R}(A)>0$. If
$|A|=1$ then ${\cal R}(A)=1/2$. Similar constraints hold for $D$. Finally
$|B||C|=1$.}

{\bf Proof}\quad From the equations satisfied by the entries we have
$A=1-D^{-1}C^{-1}DC$ and \hfil\break $D^{-1}C^{-1}DC\in S^3-\{1\}$. So $|A|\le
2$ and ${\cal R}(A)>0$. The case $|A|= 2$ cannot occur since then $A=2$ which
is the commutative case. Suppose $A=r+P$ where $r={\cal R}(A)$. Then if
$|A|^2=1=r^2+|P|^2$, we use the fact that $|1-A|^2=1=(1-r)^2+|P|^2$ to see
that ${\cal R}(A)=1/2$. Similar constraints hold for $D=1-A^{-1}B^{-1}AB$.
Finally note that $1-A=A^{-1}BAC=D^{-1}C^{-1}DC$. So $|B||C|=1$.\qed

We now consider all the switches over special subrings of $\H$ starting with
the ring, $R$, of quaternions with integer coefficients, namely
$$R=\{a+bi+cj+dk|a,b,c,d\in\Z\}.$$
By the above $A=1\pm i,\ 1\pm j,\ 1\pm k$. Note that say $A=1+i+j$ cannot
occur because $|A-1|=1$. Similarly for $D$.

The entries $B, C$ can only take the values $\pm i,\ \pm j,\ \pm k$.

Call a switch of {\bf Budapest} type if 
$$S=\pmatrix{
1+U&-V\cr V&1+U\cr}$$
where $U,V\in S^2$ and $U\perp V$.

\lemma{All switches with entries in the ring $R$ of quaternions with integer
coefficients are of Budapest type with $U,V$ lying in the set $\{\pm i,\ \pm
j,\ \pm k\}$.} \qed

Let $\xi=1/2+1/2i+1/2j+1/2k$. The {\bf Hurwitz}-ring, $H$, is defined to be
the set $$H=\{n\xi+mi+pj+rk|n,m,p,r\in\Z\}.$$ Let us extend the coefficients
of the switches to include elements of $H$. The new possible values for $A$
and $D$ are $1/2\pm i/2\pm j/2 \pm k/2$. The new possible values for $B$
and $C$ are $\pm 1/2\pm i/2\pm j/2 \pm k/2$. A computer search for all
solutions with entries in $H$ is given in table 2 at the end of the paper.

\section{Some Calculations of the Polynomials}

We now calculate the polynomials for various choices of the matrix $S$. This
not only allows us to distinguish various virtual knots but also helps to
distinguish the construction of the polynomial by $S$.

We now pick 4 switches from the more than 150 discovered by computer search.
We have been able to reduce the number of switches by 50\% from earlier
versions using the results of \cite{BF}.

$$\vbox{ \offinterlineskip \halign{\strut
\vrule \hfil \quad $#$ \ \hfil \vrule & 
       \hfil \ $\ninepoint #$ \ \hfil \vrule & 
       \hfil \ $\ninepoint #$ \ \hfil \vrule & 
       \hfil \ $\ninepoint #$ \ \hfil \vrule & 
       \hfil \ $\ninepoint #$ \ \hfil \vrule\cr
\noalign{\hrule} 
& A & B & C & D\cr 
\noalign {\hrule} 
1 & 1+i & -j & j & 1+i \cr 
\noalign {\hrule} 
2 & 1+i & -1/2+i/2+j/2-k/2  & -1/2+i/2-j/2+k/2 & 1/2+i/2+j/2+k/2 \cr 
\noalign {\hrule} 
3 & 1/2+i/2+j/2+k/2 &  -1/2+i/2-j/2+k/2 & -1/2+i/2+j/2-k/2 & 1+i \cr 
\noalign {\hrule} 
4 & 1+i & -1+i-k & -1/3+i/3+k/3 & 1/3+i/3+j2/3 \cr 
\noalign {\hrule} 
}}$$
Note that 1 is the Budapest switch, 2, and 3 have coefficients in the
Hurwitz ring and 4 is an anomalous type thrown up in the computer search.

The table below gives the corresponding normalised polynomials, $\Delta_0$,
for the virtual trefoil.
$$\vbox{ \offinterlineskip \halign{\strut
\vrule \hfil \quad $#$ \quad \hfil \vrule & 
       \hfil \ $#$ \ \hfil \vrule\cr 
\noalign{\hrule} 
Switch & \Delta_0 \cr 
\noalign {\hrule} 
1 & t^4+2t^2+1\cr 
\noalign {\hrule} 
2 & 3/4t^4+3/2t^3+9/4t^2+3/2t+3/4 \cr 
\noalign {\hrule} 
3 & 3/64t^4+3/16t^3+9/16t^2+3/4t+3/4\cr 
\noalign {\hrule} 
4 & 9t^4+12t^3+10t^2+4t+1\cr 
\noalign {\hrule} 
}}$$
For the Kishino knots, in all cases $\Delta_0$ is zero. The next table shows
$\Delta_1$, the 1-st ideal polynomial for these knots as determined by the
same set of switches.
$$\vbox{ \offinterlineskip \halign{\strut
\vrule \hfil \quad $#$ \quad \hfil \vrule & 
       \hfil \ $#$ \ \hfil \vrule & 
       \hfil \ $#$ \ \hfil \vrule & 
       \hfil \ $#$ \ \hfil \vrule\cr 
\noalign{\hrule} 
Switch & K1 & K2 & K3 \cr 
\noalign {\hrule} 
1 & t^4+5/2t^2+1 & t^4+5/2t^2+1 & t^4+5/2t^2+1 \cr 
\noalign {\hrule} 
2 & 1/2t^4+3/2t^2+1 & 2t^4+3t^2+1 & 75/64 \cr 
\noalign {\hrule} 
3 & 1/8t^4+3/4t^2+1 & 1/32t^4+3/8t^2+1 & 75/64 \cr 
\noalign {\hrule} 
4 & 3t^4+7/2t^2+1 & 27t^4+21/2t^2+1 & 58381/36450 \cr 
\noalign {\hrule}}}$$
Now consider the following knot which has trivial Jones' polynomial and
trivial fundamental rack.
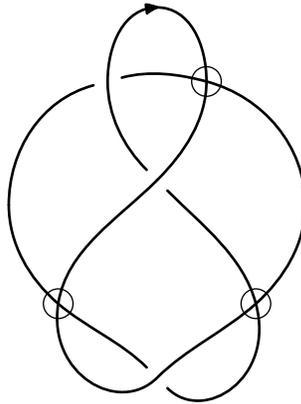
\diagram
In this case for the commutative Alexander switch
$$\Delta_0=(s-1)(st-1)(t^2-1)$$
up to multiplication by a unit.

For the quaternion switches we get
$$\vbox{ \offinterlineskip \halign{\strut
\vrule \hfil \quad $#$ \quad \hfil \vrule & 
       \hfil \ $#$ \ \hfil \vrule\cr 
\noalign{\hrule} 
Switch & \Delta_0 \cr 
\noalign {\hrule} 
1 & t^8-2t^4+1 \cr 
\noalign {\hrule} 
2 & 2t^8+4t^7+4t^6-3t^4-2t^3+t^2+2t+1 \cr 
\noalign {\hrule} 
3 & 1/512t^8+1/128t^7+1/128t^6-1/32t^5-3/32t^4+1/2t^2+t+1 \cr 
\noalign {\hrule} 
4 & 243t^8+324t^7+216t^6+36t^5-24t^4-12t^3+4t^2+4t+1 \cr 
\noalign {\hrule} 
}}$$
\np

\section{Questions and further developements}

It seems likely that the methods described above will not yield new invariants
for real knots and links. However it is quite possible that in identifying the
old invariants new properties of these will be found. This awaits a further
paper.

On the other hand these methods provide a rich tableaux of virtual knot and
link invariants. Much more can be done to understand the set of switches
described in section 8. For example are there 2$\times$2 matrices in which the
entries are polynomials in a real variable with quaternionic coefficients
other than those described in the dodge which replaces $S$ by $S(t)$? The
answer must surely be yes. However this has so far proved too much for a
computer search to cope with. We can consider the set of switches with
quaternion entries as points of some quaternionic variety. Since each point
gives rise to a polynomial and there are many such polynomials it is likely
that this variety has many components. Each component contains a copy of the
real line. But is that all? In general we have a universal non-commutative
algebra. This is given by the four relations of section 2. What can we say
about this algebra? Certainly more than is indicated by its representation
onto the quaternionic variety.

\section{ References}

\refe

\cite{A} J. W. Alexander, Topological Invariants of Knots and Links, Trans.
Amer. Math. Soc. 30 275-306 (1928)

\cite{As} Helmer Aslaksen, Quaternionic Determinants, Math. Intel. Vol 18 no.
3 (1996)

\cite{JB} J. Birman, Braids, Links and Mapping Class Groups, Princeton
University Press (1975)

\cite{BF} S. Budden and R. Fenn, The Equation $[B,(A-1)(A,B)]=0$ and
Virtual Knots and Links, preprint.

\cite{CF} R.H. Crowell and R.H. Fox, Introduction to Knot Theory. Ginn and Co.
(1963)

\cite{Dr} V. Drinfeld, On some Unsolved Problems in Quantum Group Theory,
Quantum Groups, Lectures Notes in Maths. 1510, Springer 1-8 (1990)

\cite{FJK} R. Fenn, M. Jordan, L. Kauffman,  The Birack; an
Invariant of Virtual Knots and Links, to appear in Topology and its
Applications. Available from \nl
www.maths.sussex.ac.uk///Staff/RAF/Maths/

\cite{FR} R. Fenn, C. Rourke. Racks and Links in Codimension Two. JKTR, No. 4,
343-406 (1992).

\cite{FRR} R. Fenn, R. Rimanyi and C. Rourke The braid-permutation group.
Topology 36, No.1, 123-135 (1997).

\cite{KK} N. Kamada and S. Kamada, Abstract Link Diagrams and Virtual
Knots,\nl preprint.

\cite{K} L.Kauffman. Virtual Knot Theory, European J. Comb. Vol 20, 
663-690, (1999)

\cite{M} V. O. Manturov, On Invariants of Virtual Links, Acta Math. 00 1-15
(2002)

\cite{KS} Toshimasa Kishino and Shin Satoh, A note on non-classical virtual
knots, preprint

\cite{SW}D.S. Silver and S.G. Williams, Polynomial Invariants of Virtual
Links, JKTR 12 (2003) 987-1000.

\section{Tables}

In this section we give some examples of switches found by a combination of
theory and computer search. If the switch $S$ appears then the alternatives,
$S^{-1}, S^\dagger, S^*$, are not included: neither are those obtained by a
permutation of $i, j, k$ up to sign. This pruning considerably reduces the
number of switches. The first table contains switches with coefficients in the
Hurwitz ring. Table 2 contains switches where at least two of the entries have
integer coefficients. Since $|B||C|=1$ we have restricted ourselves to integer
coefficients for $B$ and $C$ in the set $\{-1,0,1\}$. The tables have
overlaps. For example switches 1 are both variants of the Budapest switch.
The software to search for switches and calculate polynomials can be
downloaded from http://www.layer8.co.uk/maths/braids/

\newcount\tablelinenum \tablelinenum=0
\def\nextline{\global\advance\tablelinenum by 1 
              \the\tablelinenum 
             }

\bigskip
\centerline{\bf Table 1. Hurwitz Coefficients}
\bigskip \tablelinenum=0
\hbox{\kern-3em$$\vbox{ \offinterlineskip \halign{\strut
\vrule \hfil \quad $\ninepoint #$ \ \hfil \vrule & 
       \hfil \ $\ninepoint #$ \ \hfil \vrule & 
       \hfil \ $\ninepoint #$ \ \hfil \vrule & 
       \hfil \ $\ninepoint #$ \ \hfil \vrule & 
       \hfil \ $\ninepoint #$ \ \hfil \vrule\cr
\noalign{\hrule} 
& A & B & C & D\cr 
\noalign{\hrule}
\nextline & 1+i & j & -j & 1+i \cr
\noalign{\hrule}
\nextline & 1+i & 1/2-1/2i+1/2j+1/2k & 1/2-1/2i-1/2j-1/2k & 1/2+1/2i+1/2j-1/2k \cr
\noalign{\hrule}
\nextline & 1+i & -1/2+1/2i+1/2j+1/2k & -1/2+1/2i-1/2j-1/2k & 1/2+1/2i-1/2j+1/2k \cr
\noalign{\hrule}
\nextline & 1/2+1/2i+1/2j+1/2k & 1/2+1/2i-1/2j-1/2k & 1/2-1/2i+1/2j-1/2k & 1+k \cr
\noalign{\hrule}
\nextline & 1/2+1/2i+1/2j+1/2k & -1/2+1/2i+1/2j-1/2k & -1/2-1/2i+1/2j+1/2k & 1+j \cr
\noalign{\hrule}}}$$}


\centerline{\bf Table 2. At least Two Integer Coefficients}
\bigskip \tablelinenum=0
\hbox{\kern-3em$$\vbox{ \offinterlineskip \halign{\strut
\vrule \hfil \quad $\ninepoint #$ \ \hfil \vrule & 
       \hfil \ $\ninepoint #$ \ \hfil \vrule & 
       \hfil \ $\ninepoint #$ \ \hfil \vrule & 
       \hfil \ $\ninepoint #$ \ \hfil \vrule & 
       \hfil \ $\ninepoint #$ \ \hfil \vrule\cr
\noalign{\hrule} 
& A & B & C & D\cr 
\noalign{\hrule}
\nextline & 1-j & -k & k & 1-j \cr
\noalign{\hrule}
\nextline & 1+i & 1/2j+1/2k & -j-k & 1+i \cr
\noalign{\hrule}
\nextline & 1+i & 1-i-j-k & 1/4-1/4i+1/4j+1/4k & 1/2+1/2i-1/2j+1/2k \cr
\noalign{\hrule}
\nextline & 1-i & -j-k & 1/2j+1/2k & 1-i \cr
\noalign{\hrule}
\nextline & 1+j & 1-j-k & 1/3-1/3j+1/3k & 1/3+2/3i+1/3j \cr
\noalign{\hrule}
\nextline & 1-k & -1-i-j-k & -1/4+1/4i+1/4j-1/4k & 1/2-1/2i+1/2j-1/2k \cr
\noalign{\hrule}
\nextline & 1-k & -1-j-k & -1/3+1/3j-1/3k & 1/3-2/3i-1/3k \cr
\noalign{\hrule}
\nextline & 1/2+1/2i-1/2j-1/2k & -1/4+1/4i-1/4j+1/4k & -1-i-j-k & 1-j \cr
\noalign{\hrule}
\nextline & 1/2+1/2i+1/2j-1/2k & 1/4+1/4i-1/4j+1/4k & 1-i-j-k & 1+j \cr
\noalign{\hrule}
\nextline & 1/3+2/3i-1/3j & -1/3-1/3j+1/3k & -1-j-k & 1-j \cr
\noalign{\hrule}
\nextline & 1/3-2/3i+1/3k & 1/3+1/3j-1/3k & 1-j-k & 1+k \cr
\noalign {\hrule}}}$$}
\bye